\documentclass[11pt]{article}

\usepackage{amsthm, amsmath}
\usepackage{amssymb}
\usepackage[arrow,matrix,cmtip]{xy}

\theoremstyle{plain}
\newtheorem{theorem}{Theorem}[section]

\newtheorem{example}[theorem]{Example}

\theoremstyle{definition}
\newtheorem{remark}[theorem]{Remark}

\newcommand{\im}{{\rm{im}\:}}
\newcommand{\Tor}{{\rm{Tor}}}
\newcommand{\Ext}{{\rm{Ext}}}
\newcommand{\Cen}{{\rm{Cen}}}

\newcommand{\env}{{{R^{\rm{env}}}}}

\begin{document}



\title{Effective computation of $\Tor_k(M,N)$}
\author{Socorro 
Garc\'{\i}a Rom\'an
\thanks{Supported by DGUI, Consejer\'{\i}a de Educaci\'on - Gobierno de Canarias}
\and
Manuel
Garc\'{\i}a Rom\'an
}



\maketitle

\begin{abstract}
An effective method to compute a presentation of $\Tor _k (M,N)$ for modules on a not necessarily commutative algebra is proposed.
\end{abstract}



\section*{Introduction}

Unlike $\Ext_k(M,N)$ (see~\cite{smith} for an algorithm when the base ring $R$ is a
commutative algebra and~\cite{BCGL01} when $R$ is a PBW algebra), there are not known
effective methods to compute $\Tor_k(M,N)$ for a pair of $R$-modules $M$ and $N$.

In this note,
we propose an algorithm to compute a presentation of $\Tor_k(M,N)$ when $R$ is
a Poincar\'e-Birkhoff-Witt algebra (PBW algebra for short;
see~\cite[et.al.]{BCGL01, libro} for a definition and examples, including the Weyl
algebras, universal enveloping algebras of Lie algebras, the quantum plane, algebras of
quantum matrices and other iterated Ore extensions, etc.).
Since in general $\Tor_k(M,N)$ is just an abelian group when $M$
and $N$ are left $R$-modules, we ask for a two-sided structure on $M$. When $M$ is an
$R$-bimodule then $\Tor_k(M,N)$ is a left $R$-module. We show
that if, in addition, $M$ is finitely generated and $N$ is a finitely generated left
$R$-module, then effective techniques involving G\"robner bases may be used to compute
$\Tor _k(M,N)$ for any $k\ge 0$.

Besides 
the standard algorithm to compute (left) syzygies, the ingredients of our method, are
\begin{itemize}
\item an algorithm for computing a free resolution of the left
$R$-module $N$ given a finite system of generators of $N$. Such an algorithm may be
found in~\cite{libro};
\item a finite presentation of the $R$-bimodule $M$. Using the syzygy bimodule
described in~\cite{gbsb}, we show below a method for computing
such a presentation when $M$ is a centralizing $R$-bimodule and a finite system
of two-sided generators of $M$ is given;
\item theorem~\ref{computeTor};
\end{itemize}
Our algorithm, described in detail in~\ref{observaciongeneral}, follows the lines of~\cite{BCGL01}.

\section{Preliminaries}

This section is devoted to fixing notation and describing some
isomorphisms which will be used later.

Let $M$ be an $R$-bimodule.

For all $s \geq 1$, the map
$\alpha: M^s  \longrightarrow  M \otimes_R R^s ; \
\alpha(f_1,\ldots ,f_s) = \sum_{i=1}^s f_i \otimes e_i $,
where $\{e_1,..,e_s\}$ is the canonical basis of $R^s$,
is an isomorphism of left $R$-modules.

Indeed, $\alpha$ is the composition of the isomorphisms
\begin{equation}\label{isoalpha}
\begin{array}{rcccl}
M^s  & 
{\longrightarrow} &
(M \otimes_R R)^s &
{\longrightarrow} & M \otimes_R R^s \\
(f_1,\ldots ,f_s) & \mapsto & (f_1 \otimes 1,\ldots ,f_s \otimes 1) &
\mapsto & \sum_{i=1}^s f_i \otimes e_i.
\end{array}
\end{equation}

On the other hand, if $A$ is a subbimodule of $R^m$ and $B$ is a left
submodule of $R^s$, then
\begin{equation}\label{isobeta}
\begin{array}{rcl}
\beta: (R^m \otimes_R R^s) / T &
\longrightarrow & (R^m / A) \otimes_R (R^s / B) \\
(f \otimes g) +T  & \mapsto & (f +A) \otimes (g+B),
\end{array}
\end{equation}
where $T= R^m \otimes B\: +\: A \otimes_R R^s$, is an isomorphism of
left $R$-modules.

Furthermore, if $\{a_1,\ldots ,a_r\}$ is a generator system of $A$
such that $A=Ra_1+\cdots + Ra_r=a_1R+\cdots + a_rR$ and
$\{b_1,\ldots ,b_t\}$ is a generator system of $B$ as a left
$R$-module, then
\[
\{a_i \otimes b_j\: /\ 1 \leq i \leq r,\: 1 \leq j \leq t\}
\]
is a generator system of $A \otimes_R B$ as a left $R$-module, since
for all $a=\sum_{l=l}^r p_la_l \in A$ and
$b=\sum_{j=1}^t p'_jb_j \in B$ with $p_l,p'_j \in R$,
\begin{eqnarray*}
a \otimes b & = & \sum_{j,l} p_la_l \otimes p'_jb_j \\
 & = & \sum_{j,l} p_l(a_lp'_j) \otimes b_j \\
 & = & \sum_{j,l} p_l(\sum_{i=1}^r p^{jl}_ia_i) \otimes b_j
 \\
 & = & \sum_{i,j,l} p_lp^{jl}a_i \otimes b_j.
\end{eqnarray*}

>From here on, let  $0\rightarrow L \rightarrow R^m \stackrel{p_M}{\rightarrow}
M \rightarrow 0$ be a finite presentation of the $R$-bimodule $M$, and let
$H=\{h_1,\ldots ,h_r\} \subseteq R^m$ be a two-sided generator system of $L$
in such a way that $L=RH=HR$ (e.g., $H$ is a two-sided Gr\"obner basis of $L$).

\begin{remark}\label{bimodulocentralizante}
If $M$ is a centralizing bimodule, such a presentation may be computed as follows.
Recall that the $R$-bimodule $M$ is said to be centralizing
if $M$ is generated as a left $R$-module (or
equivalently, as a right $R$-module) by its centralizer
\[
\Cen_R M=\{m \in M;\, rm=mr, \forall\: r \in R\}.
\]
Let $R$ be a PBW algebra and $M$ an $R$-bimodule.
In~\cite{gbsb} we define the syzygy bimodule of $\{m_1,\ldots ,m_s\}\subseteq M$
as the kernel of the morphism of left $\env $-modules $(\env)^s \longrightarrow M;\ \ e_i \longmapsto m_i$,
where $\env$ is the enveloping algebra $R\otimes R^{\rm{op}}$. We also provide an
algorithm to compute the syzygy bimodule of $\{m_1,\ldots ,m_s\}$ .

If the $R$-bimodule $M$ is centralizing, say generated by
$\{m_1,\ldots ,m_s\}\in \Cen _R M$, then the syzygy bimodule is
$(\mathfrak{m}^s)^{-1} L$, where $\mathfrak{m}^s:(\env )^s \longrightarrow R^s$ is
given by $\mathfrak{m}(f\otimes g)= fg$ and $L$ is the kernel of the epimorphism
$p_M : R^s \longrightarrow M$. In that case, if $G$ is a (left) Gr\"obner
basis of the syzygy bimodule, then
$H= \mathfrak{m}^s (G) \setminus \{ 0 \}$ is a two-sided
Gr\"obner basis of $L$ (cf.~\cite[thm.\ 7]{gbsb}).



\end{remark}

If $N$ is a left $R$-module
and
$ 0 \rightarrow B
{\rightarrow}  R^s
{\rightarrow}
 N
\rightarrow 0$
is a presentation of $N$, then the map
\begin{equation} \label{isogamma}
\begin{array}{rccl} \gamma: & R^{ms} / \alpha^{-1}(T) & \longrightarrow & M \otimes_R N \\ & e_{ij}+ \alpha^{-1}(T) & \longmapsto & p_M\:(e'_j) \otimes p_N(e_i)
\end{array}
\end{equation}
is an isomorphism of left $R$-modules, where
$T=R^m \otimes_R B \: +\:  L \otimes_R R^s$, and
$\{e_1,\ldots ,e_{s}\}$, resp.\ $\{e'_1,\ldots ,e'_{m}\}$, resp.\
$\{e_{11},\ldots ,e_{1m},\ldots ,e_{s1},\ldots ,$ $e_{sm}\}$ are the canonical bases
of $R^s$, resp.\ $R^m$, resp.\ $R^{ms}=(R^m)^s$, being
$e_{ij}=(0_{R^m},\ldots ,\stackrel{i}{e'_j},\ldots ,0_{R^m})$.

Indeed, $\gamma$ is the composition of isomorphisms
\[
 R^{ms} / \alpha^{-1}(T) \stackrel{\bar{\alpha}}{\longrightarrow}
(R^m \otimes_R R^s) / T \stackrel{\beta}{\longrightarrow}
(R^m / L) \otimes_R (R^s / B)  \stackrel{\bar{p}_M \otimes \bar{p}_N}{\longrightarrow}
 M \otimes_R N 
\]
where $\bar{\alpha}$, resp.\ $\bar{p_M}$, resp.\ $\bar{p_N}$ is obtained by factoring
$\alpha : R^{ms} \longrightarrow (R^m \otimes_R R^s) / T$, resp.\
$p_M: R^m \longrightarrow M$, resp.\ $p_N: R^s \longrightarrow N$ through
the quotient.

Moreover, if $\{g_1,\ldots ,g_t\}$ is a generator system of $B$ as a
left $R$-module, with $g_k=(g_{k1},..,g_{ks})$, then
\begin{equation} \label{sistgen2}
\{\ (e'_jg_{k1},\ldots ,e'_jg_{ks})\ \}_{1 \leq k \leq t \atop 1 \leq j
\leq m}\ \cup\ \{\ (0_{R^m},\ldots ,\stackrel{i}{h_l},\ldots ,0_{R^m})\ \}_{1 \leq i \leq
s,\atop 1 \leq l\leq r }
\end{equation}
is a generator system of
$\alpha^{-1}(T)= \alpha^{-1}(R^m \otimes_R B) + L^s$ as a left
$R$-module,
since
\[
\{\ e'_j \otimes g_k\ \}_{ 1 \leq j \leq m,\atop 1 \leq k \leq t } \
\cup\ \{\ h_l \otimes e_i\ \}_{ 1 \leq l\leq r,\atop 1 \leq i \leq s }
\]
is a generator system of $T=R^m \otimes_R B + L \otimes_R R^s$ and
\begin{eqnarray*}
\alpha^{-1}(T) & = & \alpha^{-1}(R^m \otimes_R B)+ \alpha^{-1}(L \otimes_R R^s) \\
 & = & \alpha^{-1}(R^m \otimes_R B) + L^s \\
 & = & _R \langle\ \alpha^{-1}(e'_j \otimes g_k),\: \alpha^{-1}(h_l \otimes e_i)
\ \rangle_{i, j, k, l} \\
 & = & \: _R \langle\ (g_{k1}e'_j,\ldots ,g_{ks}e'_j),(0_{R^m},\ldots ,\stackrel{i}{h_l},\ldots ,0_{R^m})
\ \rangle_{i, j, k, l}.
\end{eqnarray*}

\section{Computing $\Tor_k(M,N)$}

Now, let $R=k\{x_1,\ldots ,x_n;\ Q, Q', \preceq\}$ be a left PBW ring and
$N$ a finitely generated left $R$-module. Let
\begin{equation}\label{resolucionN}
 \cdots
\stackrel{\partial_{k+1\ }}{\longrightarrow} R^{s_{k}}
\stackrel{\partial_{k}}{\longrightarrow} R^{s_{k-1}}
\longrightarrow \cdots \stackrel{\partial_1}{\longrightarrow}
R^{s_0} \stackrel{\partial_0}{\longrightarrow} N \longrightarrow 0
\end{equation}
be a free resolution of $N$, where $\partial _k$ is the matrix
\[ \left( \begin{array}{c} g_1^k \\ \vdots \\ g_{s_k}^k \end{array} \right)
 \]
with $g_i^k=(g_{i1}^k,\ldots ,g_{is_{k-1}}^k) \in R^{s_{k-1}}$ for $k\geq 1$
and $1\le i \le s_k$ (an algorithm to compute it may be found
in~\cite[chapter 6]{libro}).

In this section we show a method for computing ${\rm{Tor}}_k(M,N)$ when $M$ is a  finitely presented $R$-bimodule, i.e., we show how to compute the $k$-th homology
module of the complex
\begin{equation}\label{resolucionMotimes}
\cdots
\stackrel{1_M \otimes \partial_{k+1\ }}{\longrightarrow} M
\otimes_R R^{s_{k}} \stackrel{1_M \otimes
\partial_{k\ }}{\longrightarrow}
\cdots \stackrel{1_M \otimes \partial_{1\ }}
{\longrightarrow} M \otimes_R R^{s_0}
\longrightarrow 0 .
\end{equation}

\begin{remark}\label{observacionisoTor}
The case $k=0$ may be treated apart, since $\Tor_0(M, N) \cong M \otimes_R N$. 

We start with the presentations
$0\rightarrow L \rightarrow  R^m {\rightarrow} M \rightarrow 0$ and
$0\rightarrow \ker{\partial_0} \rightarrow  R^{s_0} {\rightarrow} N
\rightarrow 0$ of the $R$-bimodule $M$, resp.\ of the left
$R$-module $N$.

Since  $ \{ g_1^1,\ldots ,g_{s_1}^1 \}$ is a generator system
of $\ker{\partial_0}=\im{\partial_1}$ as a left $R$-module and
\[
\gamma: R^{ms_0} / (\alpha^{-1}(R^m \otimes_R
\ker{\partial_0})+L^{s_0}) \longrightarrow \Tor_0( M , N),
\]
as in~(\ref{isogamma}) is an isomorphism of left $R$-modules, we completely describe
a presentation of $\Tor_0( M , N)$ by giving the generator system
\[
\{\ (e'_jg_{l1}^1,\ldots ,e'_jg_{ls_0}^1)\ \}_ {1 \leq l \leq s_1,\atop 1 \leq j
\leq m}\ \cup\ \{\ (0_{R^m},\ldots ,\stackrel{i}{h_l},\ldots ,0_{R^m})\ \}_{ 1 \leq i
\leq s_0,\atop 1 \leq l\leq r}
\]
of $\alpha^{-1}(R^m \otimes_R \ker{\partial_0})+L^{s_0}$, where
$\{ h_1,\ldots ,h_r\}$ is a two-sided Gr\"obner basis of $M$
(see~(\ref{sistgen2})).
\end{remark}

Let us return to the general case. Consider again the presentation
$0\rightarrow L \rightarrow  R^m \stackrel{p_M}{\rightarrow} M \rightarrow 0$ of the
$R$-bimodule $M$ and the free resolution (\ref{resolucionN}) of $N$.

For all $k\ge 1$, let $d_k=\gamma_{k-1}^{-1}\circ (1_M \otimes \partial_k)\circ
\gamma_k$, where $\gamma_k: R^{ms_k} /L^{s_k}  \longrightarrow  M \otimes_R R^{s_k}$
is the isomorphism defined as in~(\ref{isogamma}) with $N=R^{s_k}$.

Clearly, the complex
\begin{equation}\label{complejodk}
\cdots \longrightarrow  R^{ms_k} / L^{s_k} \stackrel{d_k}{\longrightarrow}
R^{ms_{k-1}} / L^{s_{k-1}} \stackrel{d_{k-1}}{\longrightarrow} \cdots
\stackrel{d_1}{\longrightarrow} R^{ms_0} / L^{s_0}  \longrightarrow  0
\end{equation}
is isomorphic to the one in (\ref{resolucionMotimes}), so
$\Tor_{k}(M,N)$ may be computed as the homology of~(\ref{complejodk}).

By definition, for all $1 \leq i \leq s_k$, $1 \leq j \leq m$,
\begin{equation}\label{homomorfismodk}
\begin{array}{rcl}
d_k(e_{ij}+L^{s_k}) & = &
\gamma_{k-1}^{-1}(1_M \otimes \partial_k)\: \gamma_k(e_{ij}+L^{s_k})  \\
& = & \gamma_{k-1}^{-1}((1_M \otimes \partial_k)\: (p_M(e'_j) \otimes e_i))  \\
& = & \gamma_{k-1}^{-1}(p_M(e'_j) \otimes (g_{i1}^k,\ldots ,g_{is_{k-1}}^k)) \\
& = & (g_{i1}^ke'_j,\ldots ,g_{is_{k-1}}^ke'_j)+L^{s_{k-1}}.   \\
\end{array}
\end{equation}
Let $\widetilde{d_k} : R^{ms_k} \longrightarrow R^{ms_{k-1}}$
be the block-built matrix
\[
A_k = \left( \begin{array}{c|c|c}
g_{11}^kI_m & \cdots & g_{1s_{k-1}}^kI_m \\
\hline \vdots & & \vdots \\
\hline g_{s_k1}^kI_m & \cdots & g_{s_ks_{k-1}}^kI_m
\end{array}\right) \in {\cal{M}}_{ms_k \times ms_{k-1}}(R),
\]
where $I_m$ denotes the $m \times m$-identity matrix. Since $A_k$ is
built of blocks which are elements of $R$ times the identity matrix, we have $\widetilde{d_k}(L^{s_k})\subseteq L^{s_{k-1}}$, and the diagram
\[\xymatrix{
R^{ms_k} \ar[r]^{\pi_k} \ar[d]_{\widetilde{d_k}} & R^{ms_k} / L^{s_k} \ar[d]^{d_k} \\
R^{ms_{k-1}} \ar[r]_{\pi_{k-1\ \ \ }} & R^{ms_{k-1}} /
L^{s_{k-1}}}
\]
is commutative.

The above discussion proves the following result:

\begin{theorem} \label{computeTor}
With the previous notation, for all $k \geq 1$
\begin{enumerate}
\item $\ker d_k=\ker \pi_{k-1}\widetilde{d_k} / L^{s_k}$;
\item $\im d_k=\im \pi_{k-1}\widetilde{d_k} \subseteq R^{ms_{k-1}} / L^{s_{k-1}}$
is generated by
\[
\{ (g_{i1}^ke'_j,\ldots ,g_{is_{k-1}}^ke'_j)+ L^{s_{k-1}}
\}_{1\leq i \leq s_k\atop 1 \leq j \leq m}
\]
as a left $R$-module (note that
$(g_{i1}^ke'_j,\ldots ,g_{is_{k-1}}^ke'_j)$ is the $(j+m(i-1))$-th row of $A_k$);
\item 
$\Tor_k(M,N)  \cong  \ker d_k / \im d_{k+1}
 = \ker \pi_{k-1}\widetilde{d_k} / (_R \langle
{\rm{rows\ of\ }} A_{k+1} \rangle + L^{s_k})$.
\end{enumerate}
\end{theorem}

\begin{remark}\label{observaciongeneral}
To compute $\Tor_k(M,N)$ we start with a finite presentation
$0\rightarrow L \rightarrow R^m \stackrel{p_M}{\rightarrow} M \rightarrow 0$
of the $R$-bimodule $M$ (say, e.g., we have computed a two-sided Gr\"obner basis
$\{h_1,\ldots ,h_r\} \subseteq R^m$ of the $R$ bimodule $L$) and a
free resolution of the left $R$-module $N$ as (\ref{resolucionN}).

The matrix $A_k$ is block-built as above, using the matrix
\[ \partial _k = \left( \begin{array}{c} g_1^k \\ \vdots \\ g_{s_k}^k \end{array} \right)
. \]

The set $\{\ (0_{R^m},\ldots ,\stackrel{i}{h_l},\ldots ,0_{R^m})\
\}_{i,l=1}^{s_k,r}$
is a generator system (in fact, it is a two-sided Gr\"obner
basis when $\{h_1,\ldots ,h_r\}$ so is) of $L^{s_k}$ as an $R$-bimodule.

Then, we compute the kernel of $ \pi_{k-1}\widetilde{d_k}: R^{ms_k}  \longrightarrow R^{ms_{k-1}} / L^{s_{k-1}}$ using syzygies.
Indeed (see \cite{BCGL01}), let $H_k$ be the matrix
\[ H_k=\left(\frac{A_k}{A'_k}\right) \in {\cal{M}}_{(ms_k+rs_{k-1})\times
{ms_{k-1}}}(R),\]
where $A'_k$ is the matrix whose rows are the generators of
$L^{s_{k-1}}$ as a left $R$-module. Then, if
\[
Syz^l(H_k)=\: _R \langle p_1,...,p_l \rangle,\ \ {\rm{\ with\ }}
p_i=(p'_i,p''_i) \in R^{ms_k} \times R^{rs_{k-1}},
\]
then $\ker \pi_{k-1}\widetilde{d_k}\ =\ _R \langle p'_1,...,p'_l \rangle$.
\end{remark}

\begin{example}
Let $R=U(\mathfrak{sl}(2))$, the
universal enveloping algebra of the Lie algebra of traceless $2
\times 2$-matrices, where $k=\mathbb{C}$ (or $k=\mathbb{Q}$). We know (see, e.g.,
\cite{libro}) that $R$ is the PBW algebra $k\{x,y,z;\: Q\, \preceq_{\omega}\}$ with
$ Q= \{\: yx-xy+z,\ zx-xz-2x,\ zy-yz+2y\: \}$ and, say, $\omega=(1,2,2)$.

Let $N = R^2 / B$, where $B$ is the left $R$-module generated
by
\[
g_1=(y^3,x),\ g_2=(y,xz),\ g_3=(0,xy^2z-2yz^2+2yz-x).
\]
The left syzygy module $Syz^l(g_1,g_2,g_3)$ is generated by
$g=(1,-y^2,1) \in R^3$, and hence
\begin{equation}\label{resNejemplo}
0  \longrightarrow R \stackrel{\delta_1}{\longrightarrow} R^3
\stackrel{\delta_0}{\longrightarrow} N \longrightarrow 0,
\end{equation}
where $\partial_0 =\left(\begin{array}{c} g_1 \\ g_2 \\ g_3 \end{array} \right)$ and
$\partial_1=\left( g \right)$, is a free resolution of $N$.

Let $L$ be the $R$-bimodule generated by $\{(C,1),(1,C)\}$, where
$C$ is the Casimir element $z^2/2+2xy-z$ (a well known central element of
$U(\mathfrak{sl}(2))$), and let $M=R^2/L$.

For all $k \geq 2$ we have $\Tor_k(M,N)=0$, as the free resolution
(\ref{resNejemplo}) of $N$ has length 2.

For $k=0$, we have (see~\ref{observacionisoTor})
$\Tor_0(M,N) \cong R^6\ /\ (\alpha^{-1}(R^2 \otimes_R B)+L^3)$, and
\begin{eqnarray*}
\lefteqn{\{\  (1,0,-y^2,0,1,0)\ ,\ (0,1,0,-y^2,0,1)\ ,\
(C,1,0,0,0,0)\ ,\ (0,0,C,1,0,0)\ ,}  \\
& & (0,0,0,0,C,1)\ ,\  (1,C,0,0,0,0)\ ,\ (0,0,1,C,0,0)\ , \ (0,0,0,0,1,C) \ \}
\end{eqnarray*}
is a generator system of $\alpha^{-1}(R^2 \otimes_R B)+L^3$ as a left $R$-module.

For $k=1$, we have $\Tor_1(M,N) \cong
\ker(\pi_{0}\widetilde{d_1})\ /\ (_R \langle {\rm{rows\ of\ }} A_{2} \rangle + L) $, but
in this particular example $A_2=0$ since $\partial_2=0$.

As pointed out in~\ref{observaciongeneral}, $\ker \pi_{0}\widetilde{d_1}$ is obtained
from the left syzygy module of the rows of
\[ A_1=\left(\begin{array}{cccccc} 1 & 0 & -y^2 & 0 & 1 & 0
\\ 0 & 1 & 0 & -y^2 & 0 & 1 \end{array}\right) \]
and the generators of $L^3$. Indeed, by picking up just the first two
components of each element of the generator
system
\[
\begin{array}{cll}
\{  (z^2+4xy-2z,2,-2,2y^2,-2,0,0,0)\ ,  &  \\
\   (1,z^2/2+2xy-z,0,0,0,-1,y^2,-1)\ ,  &   \\
    \multicolumn{2}{l}{\ (8xy-4z,-4z^4-4xyz^2+2z^3+4,-4,4y^2,-4,2y^2,} \\
    & -2y^2z^2+16y^2z-32y^2,2z^2) \}
\end{array}
\]
of the left syzygy module, we obtain the generator system
\[ \{ \ (2C,2)\ ,\ (1,C) \ , \ (8xy-4z,-z^4-4xyz^2+2z^3+4)\ \} \]
of $\ker \pi_{0}\widetilde{d_1}$.
Therefore, $\Tor_1(M,N) \cong \ker \pi_{0}\widetilde{d_1} / L = 0 $, since
$(8xy-4z,-z^4-4xyz^2+2z^3+4)\in L$ (one may check this out by dividing the element
by a two-sided Gr\"obner basis of $L$).
\end{example}

\section*{Acknowledgement}

Both authors would like to thank Pepe Bueso for his useful suggestions.


%


\end{document}